\documentclass{article}
\usepackage{amssymb}
\usepackage{graphicx}
\newtheorem{theorem}{Theorem}
\newtheorem{lemma}{Lemma}
\newtheorem{example}{Example}
\newtheorem{definition}{Definition}
\newtheorem{corollary}{Corollary}
\newenvironment{proof}[1][Proof]{\textbf{#1.} }{\ \rule{0.5em}{0.5em}}
\title{LULU operators and locally $\delta$-monotone approximations}
\author{ Roumen Anguelov\\Department of Mathematics and Applied Mathematics \\
University of Pretoria \\
Pretoria 0002, SOUTH AFRICA \\
E-mail: roumen.anguelov@up.ac.za} \date{}
\begin{document}
\maketitle

\begin{abstract}
The LULU operators, well known in the nonlinear multiresolution
analysis of sequences, are extended to functions defined on
continuous domain, namely, a real interval
$\Omega\subseteq\mathbb{R}$. Similar to their discrete
counterparts, for a given $\delta>0$ the operators $L_\delta$ and
$U_\delta$ form a fully ordered semi-group of four elements. It is
shown that the compositions $L_\delta\circ U_\delta$ and
$U_\delta\circ L_\delta$ provide locally $\delta$-monotone
approximations for the bounded real functions defined on $\Omega$.
The error of approximation is estimated in terms of the modulus of
nonmonotonicity.
\end{abstract}

\section{Introduction}

The LULU operators remove impulsive noise before a signal
extraction from a sequence. They are computationally convenient
and conceptually simpler compared to the the median smoothers
usually considered to be the "basic" smoothers. The LULU operators
have particular properties, e.g. they are fully trend preserving,
\cite{RohwerQM2004}, preserve the total variation,
\cite{RohwerQM2002}, etc., which make them an essential tool for
multiresolution analysis of sequences. Furthermore, it was
demonstrated during the last decade or so that these operators,
being specific cases of morphological filters, \cite{Serra}, have
a critical role in the analysis and comparison of nonlinear
smoothers, \cite{Rohwerbook}.

We extend the LULU theory from sequences to functions on a
continuous domain, namely, a real interval $\Omega$. The existing
LULU theory can be considered as a particular case in this new
development, since the discrete LULU operators can be equivalently
formulated for splines of order 1 or order 2 on the integer
partition of real line, the one-to-one mapping being given by the
B-spline basis.

Given a sequence $\xi=(\xi_i)_{i\in\mathbb{N}}$ and
$n\in\mathbb{N}$ the operators $L_n$ and $U_n$ are defined as
follows
\begin{eqnarray*}
(L_n
\xi)_i=\max\{\min\{\xi_{i-n},...,\xi_i\},...,\min\{\xi_{i},...,\xi_{i+n}\}\},\
i\in\mathbb{N}\\
(U_n
\xi)_i=\min\{\max\{\xi_{i-n},...,\xi_i\},...,\max\{\xi_{i},...,\xi_{i+n}\}\},\
i\in\mathbb{N}
\end{eqnarray*}
In analogy with the above discrete LULU operators, for a given
$\delta>0$ the basic smoothers $L_\delta$ and $U_\delta$ in the
LULU theory are defined for functions on $\Omega$ through the
concepts of the so called lower and upper $\delta$-envelopes of
these functions. These definitions are given in Section 2, where
it is also shown that the operators $L_\delta$ and $U_\delta$
preserve essential properties of their discrete counterparts. In
particular, the operators $L_\delta$ and $U_\delta$ generate
through composition a fully ordered four element semi-group also
called a strong LULU structure. This issue is dealt with in
Section 3. In section 4 we define the concept of local
$\delta$-monotonicity and show that the compositions
$L_\delta\circ U_\delta$ and $U_\delta\circ L_\delta$ are
smoothers in the sense that the resulting functions are locally
$\delta$-monotone. The errors of approximation of real functions
$f$ by the these compositions are estimated in terms of the
modulus of nonmonotonicity $\mu(f,\delta)$.

\section{The basic smoothers $L_\delta$ and $U_\delta$}

Let $ \mathcal{A}(\Omega )$ denote the set of all bounded real
functions defined on the real interval
$\Omega\subseteq\mathbb{R}$. Let $B_{\delta}(x)$ denote the closed
$\delta$-neighborhood of $x$ in $\Omega $, that is, $B_{\delta
}(x)=\{y\in \Omega :|x-y|\leq\delta \}$.
The pair of mappings $I$, $S:\mathcal{A}(\Omega )\rightarrow \mathcal{A}%
(\Omega )$ defined by
\begin{eqnarray}
I(f)(x) &=&\sup_{\delta >0}\inf \{f(y):y\in B_{\delta }(x)\}, \
x\in\Omega,
\label{lbf} \\
S(f)(x) &=&\inf_{\delta >0}\sup \{f(y):y\in B_{\delta }(x)\}, \
x\in\Omega. \label{ubf}
\end{eqnarray}%
are called lower Baire, and upper Baire operators, respectively, \cite%
{Sendov}. We consider on $\mathcal{A}(\Omega)$ the point-wise
defined partial order, that is, for any
$f,g\in\mathcal{A}(\Omega)$
\begin{equation}\label{forder}
f\leq g\Longleftrightarrow f(x)\leq g(x),\ x\in\Omega.
\end{equation}
Then the lower and upper Baire operators can be defined in the
following equivalent way. For every $f\in \mathcal{A}(\Omega )$
the function $I(f)$ is the maximal lower semi-continuous function
which is not greater than $f$. Hence, it is also called lower
semi-continuous envelope. In a similar way, $S(f)$ is the smallest
upper semi-continuous function which is not less than $f$ and is
called the upper semi-continuous envelope of $f$. In analogy with
$I(f)$ and $S(f)$ we call the functions
\begin{eqnarray}
I_{\delta }(f)(x) &=&\inf \{f(y):y\in B_{\delta }(x)\},\
x\in\Omega,
\label{ldenvelop} \\
S_{\delta }(f)(x) &=&\sup \{f(y):y\in B_{\delta }(x)\},\
x\in\Omega. \label{udenvelop}
\end{eqnarray}%
a lower $\delta$-envelope of $f$ and an upper $\delta $-envelope
of $f$, respectively.

It is easy to see from (\ref{ldenvelop}) and (\ref{udenvelop})
that for every $\delta_1,\delta_2>0$
\begin{equation}
I_{\delta _{1}}\circ I_{\delta _{2}}=I_{\delta _{1}+\delta _{2}}\
\ ,\ \ \  S_{\delta _{1}}\circ S_{\delta _{2}}=S_{\delta
_{1}+\delta _{2}}\label{IScomp}
\end{equation}
Furthermore, the operators $I_\delta$ and $S_\delta$, $\delta>0$,
as well as $I$ and $S$ are all monotone increasing with respect to
the order (\ref{forder}), that is, for every
$f,g\in\mathcal{A}(\Omega)$
\begin{equation}\label{IdSdmon}f\leq g\Longrightarrow I_\delta(f)\leq
I_\delta(g),\ S_\delta(f)\leq S_\delta(g),\ I(f)\leq I(g),\
S(f)\leq S(g).
\end{equation}

The following operators can be considered as continuous analogues
of the discrete $LULU$ operators given in the Introduction:
\[
L_{\delta }=S_{\frac{\delta}{2} }\circ I_{\frac{\delta}{2} }\ ,\ \
U_{\delta }=I_{\frac{\delta}{2} }\circ S_{\frac{\delta}{2}}\ .
\]
We will show that these operators have similar properties to their
discrete counterparts. Let us note that they inherit monotonicity
with respect of the functional argument from the operators
$I_\delta$ and $S_\delta$, (\ref{IdSdmon}), that is, for
$f,g\in\mathcal{A}(\Omega)$
\begin{equation}\label{LdUdfmon}f\leq g\Longrightarrow L_\delta(f)\leq
L_\delta(g),\ U_\delta(f)\leq U_\delta(g).
\end{equation}

\begin{theorem}\label{tLdfleqf}
For every $f\!\in\!\mathcal{A}(\Omega)$ and $\delta\!>\!0$ we have
$L_\delta(f)\leq f$, $U_\delta(f)\geq f$.
\end{theorem}
\begin{proof}
Let $f\in\mathcal{A}(\Omega)$, $\delta>0$. For any $x\in\Omega$ it
follows from the definition of $I_\delta$ that
$I_\frac{\delta}{2}(f)(y)\leq f(x),\ y\in B_\frac{\delta}{2}(x)$.
Therefore
$L_\delta(f)(x)=S_\frac{\delta}{2}(I_\frac{\delta}{2}(f))(x)=\sup\{I_\frac{\delta}{2}(f)(y):y\in
B_\frac{\delta}{2}(x)\}\leq f(x)$, $x\in\Omega$. The second
inequality in the theorem is proved in a similar way.
\end{proof}

\begin{theorem}\label{tLdUdmon}
The operator $L_\delta$ is monotone increasing on $\delta$ while
the operator $U_\delta$ is monotone decreasing on $\delta$, that
is, for any $f\in\mathcal{A}(\Omega)$ and $0<\delta_1\leq\delta_2$
we have $L_{\delta_1}(f)\leq L_{\delta_2}(f)$,
$U_{\delta_1}(f)\geq U_{\delta_2}(f)$.
\end{theorem}
\begin{proof}
Let $\delta_2>\delta_1>0$. Using the properties (\ref{IScomp}) the
operator $L_{\delta_2}$ can be represented in the form $
L_{\delta_2}=S_\frac{\delta_2}{2}\circ
I_\frac{\delta_2}{2}=S_\frac{\delta_1}{2}\circ
S_{\frac{\delta_2-\delta_1}{2}}\circ
I_{\frac{\delta_2-\delta_1}{2}} \circ
I_\frac{\delta_1}{2}=S_\frac{\delta_1}{2}\circ
L_{\delta_2-\delta_1} \circ I_\frac{\delta_1}{2}$. It follows from
Theorem \ref{tLdfleqf} that for every $f\in\mathcal{A}(\Omega)$ we
have $L_{\delta_2-\delta_1}(I_\frac{\delta_1}{2}(f))\leq
I_\frac{\delta_1}{2}(f)$. Hence using the monotonicity of the
operator $S_\delta$ given in (\ref{IdSdmon}) we obtain
$L_{\delta_2}(f)=S_\frac{\delta_1}{2}(L_{\delta_2-\delta_1}(I_\frac{\delta_1}{2}(f)))\leq
S_\frac{\delta_1}{2}(I_\frac{\delta_1}{2}(f))=L_{\delta_1}(f)$,
$f\in\mathcal{A}(\Omega)$. The inequality $U_{\delta_1}(f)\geq
U_{\delta_2}(f)$ is proved in a similar way.
\end{proof}

The next lemma is useful in dealing with compositions of
$I_\delta$ and $S_\delta$.
\begin{lemma}\label{tIdSdId} We have
$I_\delta\circ S_\delta\circ I_\delta=I_\delta$, $S_\delta\circ
I_\delta\circ S_\delta=S_\delta$.
\end{lemma}
\begin{proof}
Using the monotonicity of $I_\delta$, see (\ref{IdSdmon}), and
Theorem \ref{tLdfleqf} for $f\in \mathcal{A}(\Omega)$ we have
$(I_\delta\circ S_\delta\circ I_\delta)(f)=I_\delta
(L_{2\delta}(f))\leq I_\delta (f)$. On the other side, applying
Theorem \ref{tLdfleqf} to $U_{2\delta}$ we obtain $(I_\delta\circ
S_\delta\circ I_\delta)(f)=U_{2\delta}(I_\delta(f))\geq
I_\delta(f)$. Therefore $(I_\delta\circ S_\delta\circ
I_\delta)(f)=I_\delta(f)$, $f\in\mathcal{A}(\Omega)$. The second
equality is proved similarly.
\end{proof}

\begin{theorem}\label{tLdUdabsorb}
For every $\delta_1,\delta_2>0$ we have $L_{\delta_1 }\circ
L_{\delta_2 } =L_{\max\{\delta_1,\delta_2\}}$ and $U_{\delta_1
}\circ U_{\delta_2 }=U_{\max\{\delta_1,\delta_2\}}$.
\end{theorem}
\begin{proof}
We will only prove the first equality since the proof of the
second one is done in a similar manner. Let first
$\delta_2>\delta_1>0$. Using the property (\ref{IScomp}) and Lemma
\ref{tIdSdId} we obtain
\begin{eqnarray*}
L_{\delta_1}\circ L_{\delta_2}&=&(S_\frac{\delta_1}{2}\circ
I_\frac{\delta_1}{2})\circ(S_\frac{\delta_2}{2}\circ
I_\frac{\delta_2}{2})\ =\ (S_\frac{\delta_1}{2}\circ
I_\frac{\delta_1}{2}\circ S_\frac{\delta_1}{2})\circ
(S_{\frac{\delta_2-\delta_1}{2}}\circ I_\frac{\delta_2}{2})\\
&=& S_\frac{\delta_1}{2}\circ S_{\frac{\delta_2-\delta_1}{2}}\circ
I_\frac{\delta_2}{2}\ =\ S_\frac{\delta_2}{2}\circ
I_\frac{\delta_2}{2}\ =\ L_{\delta_2}.
\end{eqnarray*}
If $\delta_1>\delta_2>0$ in a similar way we have
\begin{eqnarray*}
L_{\delta_1}\circ L_{\delta_2}&=&(S_{\delta_1}\circ
I_{\delta_1})\circ(S_{\delta_2}\circ I_{\delta_2})\ =\
(S_{\delta_1}\circ I_{\delta_1-\delta_2})\circ (I_{\delta_2}\circ
S_{\delta_2}\circ I_{\delta_2})\\
&=&S_{\delta_1}\circ I_{\delta_1-\delta_2}\circ I_{\delta_2}\ =\
S_{\delta_1}\circ I_{\delta_1}\ =\ L_{\delta_1}.
\end{eqnarray*}
The proof in the case when $\delta_2\!=\!\delta_1\!>\!0$ follows
from either of the above identities where $S_{\delta_2-\delta_1}$
or $I_{\delta_1-\delta_2}$ respectively are  replaced by the
identity operator.
\end{proof}

Important properties of smoothing operators are their idempotence
and co-idempotence. Hence the significance of the next theorem.
\begin{theorem}\label{tLdUdidemp}
The operators $L_\delta$ and $U_\delta$ are both idempotent and
co-idempotent, that is, $L_{\delta }\circ L_{\delta }
=L_{\delta}$, $U_{\delta }\circ U_{\delta } =U_{\delta}$,
$(id-L_{\delta })\circ (id-L_{\delta })=id-L_{\delta }$,
$(id-U_{\delta })\circ (id-U_{\delta })=id-U_{\delta }$, where
$id$ denotes the identity operator.
\end{theorem}
\begin{proof}
The idempotence of $L_\delta$ and $U_\delta$ follows directly from
Theorem \ref{tLdUdabsorb}. The co-idempotence of the operator
$L_\delta$ is equivalent to $L_\delta\circ(id-L_\delta)=0$. Using
the first inequality in Theorem \ref{tLdfleqf} one can easily
obtain $L_\delta\circ(id-L_\delta)\geq 0$. Hence, for the
co-idempotence of $L_\delta$ it remains to show that
$L_\delta\circ(id-L_\delta)\leq 0$. Assume the opposite. Namely,
there exists a function $f\in\mathcal{A}(\Omega)$ and $x\in\Omega$
such that $(L_\delta\circ(id-L_\delta))(f)(x)>0$. Let
$\varepsilon>0$ be such that
$(L_\delta\circ(id-L_\delta))(f)(x)>\varepsilon>0$. Using the
definition of $L_\delta$ the above inequality implies that there
exists $y\in B_\frac{\delta}{2}(x)$ such that for every $z\in
B_\frac{\delta}{2}(y)$ we have $(id-L_\delta)(f)(z)>\varepsilon$,
or equivalently %
\begin{equation}\label{ineq3tidemp}
f(z)>L_\delta(f)(z)+\varepsilon,\ z\in B_\frac{\delta}{2}(y).
\end{equation}
For every $z\in B_\frac{\delta}{2}(y)$ we also have
$L_\delta(f)(z)\geq I_\frac{\delta}{2}(f)(y)=\inf\{f(t):t\in
B_\frac{\delta}{2}(y)\}$. Hence there exists $t\in
B_\frac{\delta}{2}(y)$ such that
$f(t)<I_\frac{\delta}{2}(f)(y)+\varepsilon\leq
L_\delta(f)(z)+\varepsilon$, $z\in B_\frac{\delta}{2}(y)$. Taking
$z=t$ in the above inequality we obtain
$f(t)<L_\delta(f)(t)+\varepsilon$, which contradicts
(\ref{ineq3tidemp}). The co-idempotence of $U_\delta$ is proved in
a similar way.
\end{proof}

\begin{example} The figures below illustrate graphically the smoothing
effect of the operators $L_\delta$, $U_\delta$ and their
compositions. The graph of function $f$ is given by dotted lines.

\begin{center}
\includegraphics[scale=0.25]{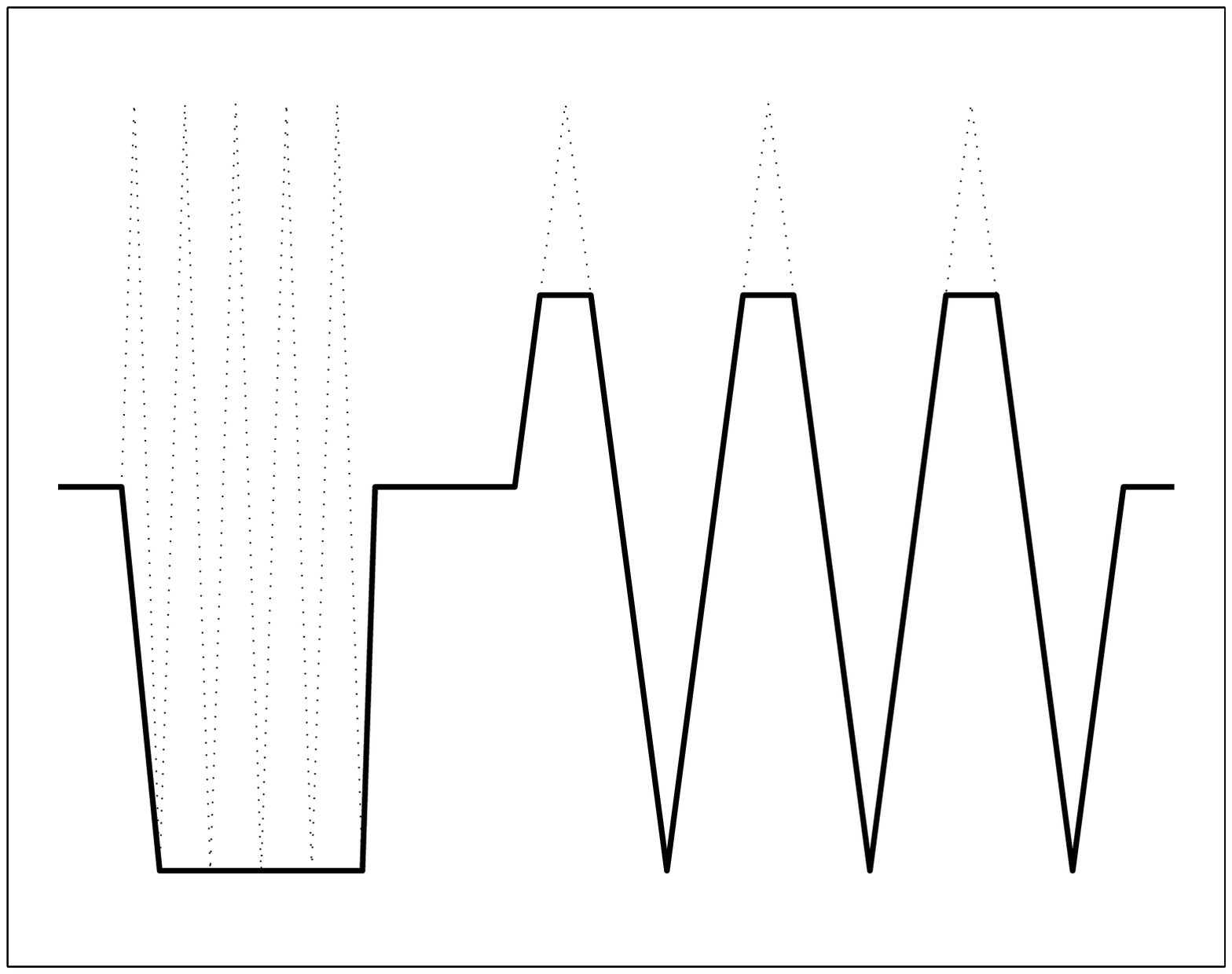}\hspace{5mm}
\includegraphics[scale=0.25]{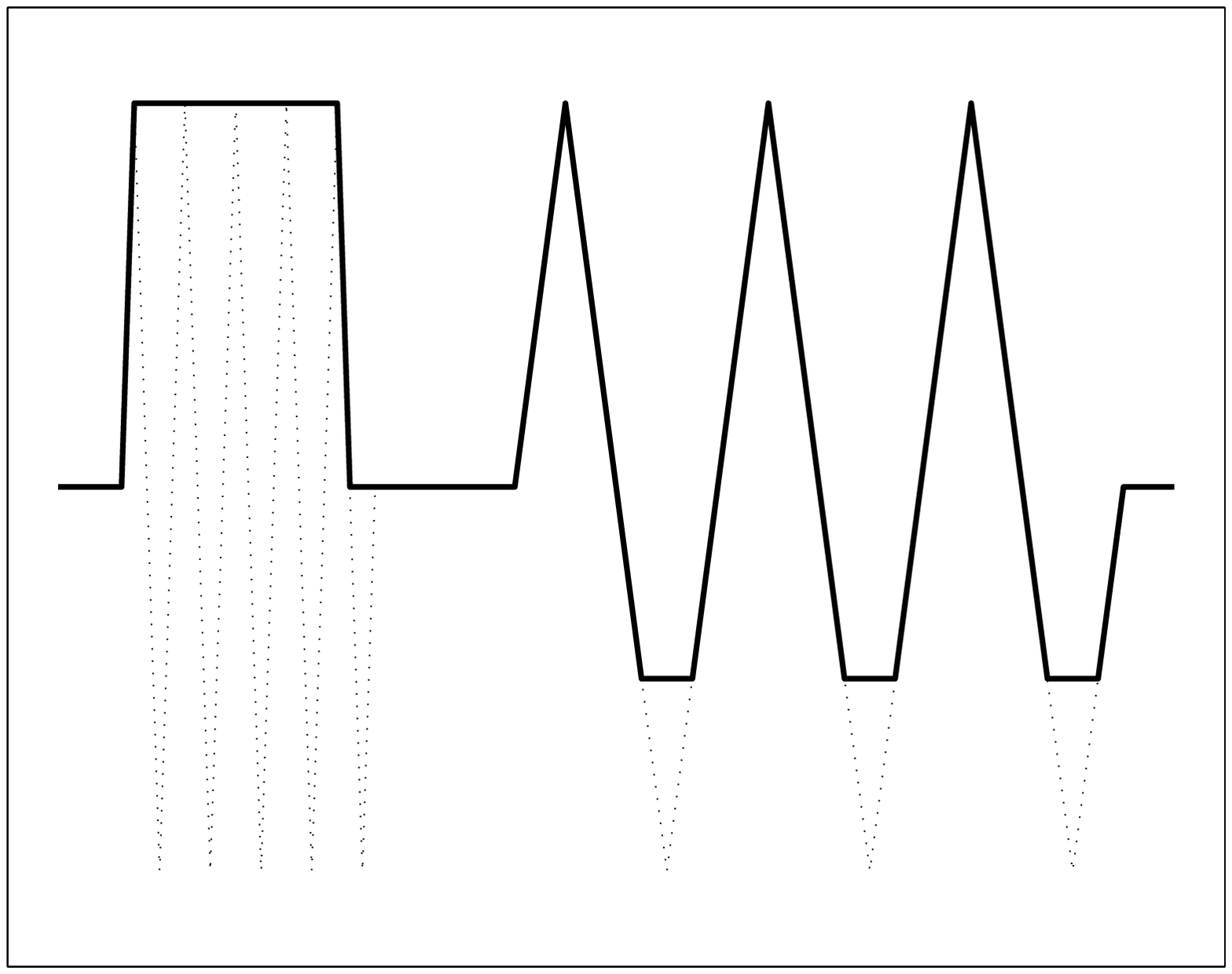}\\\nopagebreak
The functions $L_\delta(f)$ and $U_\delta(f)$\\[6pt]
\includegraphics[scale=0.25]{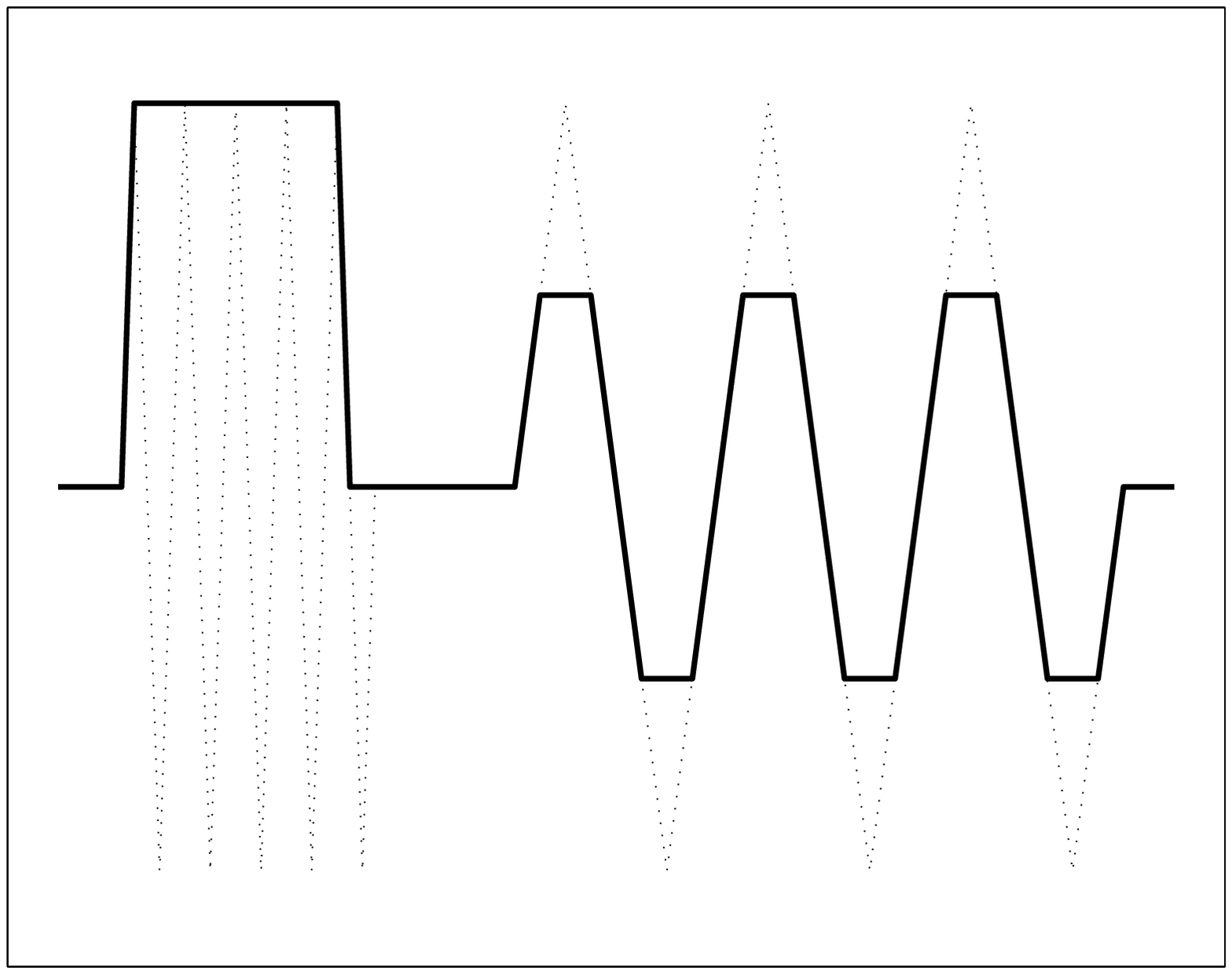}\hspace{5mm}
\includegraphics[scale=0.25]{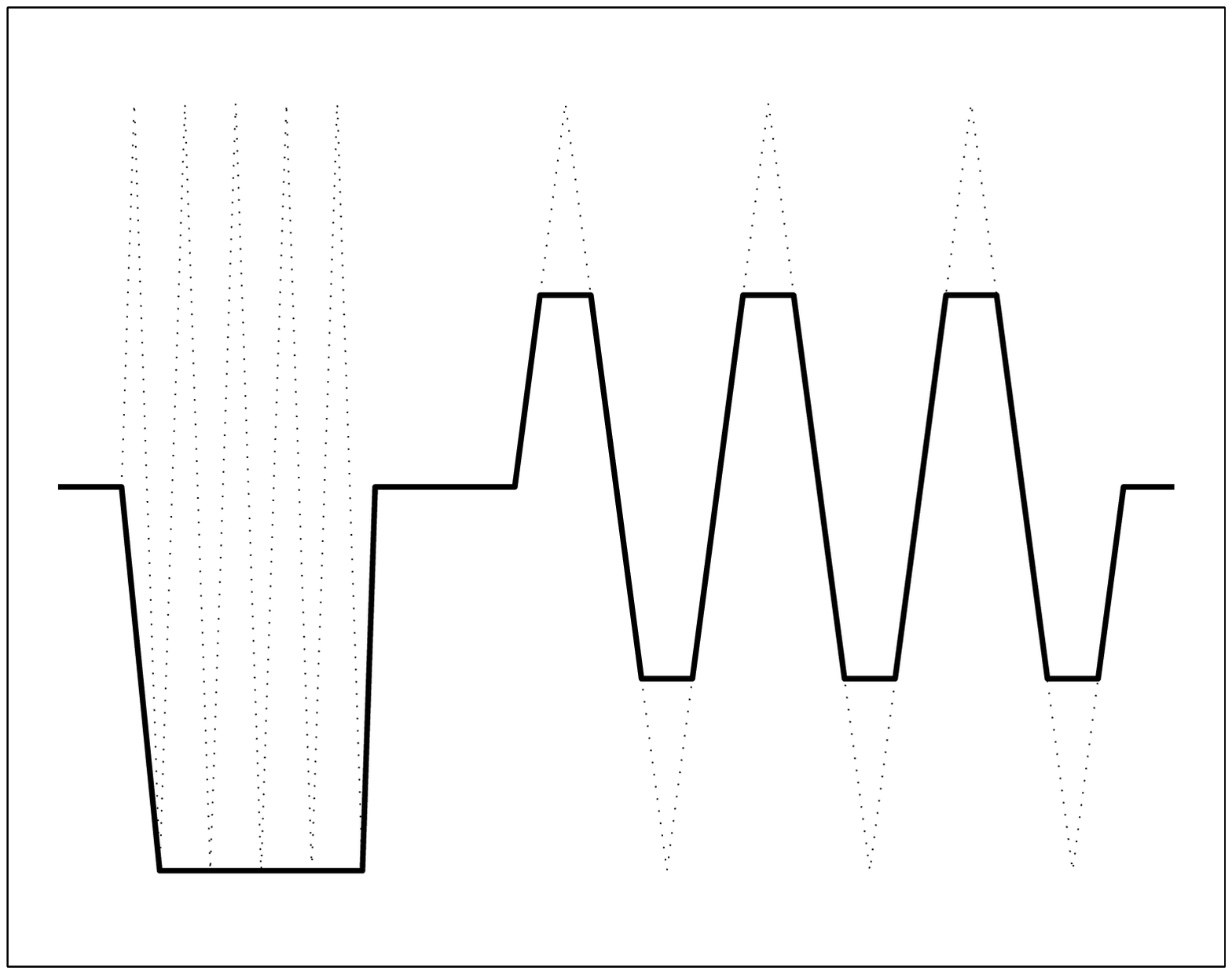}\\\nopagebreak
The functions $(L_\delta\circ U_\delta)(f)$ and $(U_\delta\circ
L_\delta)(f)$
\end{center}
\end{example}

The operator $L_\delta$ smoothes  the function $f$ from above by
removing sharp picks while the operator $U_\delta$ smoothes the
function $f$ from below by removing deep depressions. The
smoothing effect of the compositions $L_\delta\circ U_\delta$ and
$U_\delta\circ L_\delta$ can be described in terms of the local
$\delta$-monotonicity discussed in the Section 4. Note that
$L_\delta\circ U_\delta$ and $U_\delta\circ L_\delta$ resolve
ambiguities in a different way; $L_\delta\circ U_\delta$ treats
oscillations of length less then $\delta$ as upward impulses and
removes them while $U_\delta\circ L_\delta$ considers such
oscillations as downward impulses which are accordingly removed.
The inequality $(U_\delta\circ L_\delta)(f)\leq(L_\delta\circ
U_\delta)(f)$ which is observed here will be proved in the next
section for $f\in\mathcal{A}(\Omega)$.

\section{The LULU semi-group}

In this section we consider the set of the operators $L_\delta$
and $U_\delta$ and their compositions. For operators on
$\mathcal{A}(\Omega)$ we consider the point-wise defined partial
order. Namely, for operators $P$, $Q$ on $\mathcal{A}(\Omega)$ we
have
\[
P\leq Q\ \Longleftrightarrow \ P(f)\leq Q(f),\
f\in\mathcal{A}(\Omega).
\]
Then the inequalities in Theorem  \ref{tLdfleqf} can be
represented in the form
\begin{equation}\label{LleqidleqU}
L_\delta\leq id\leq U_\delta,
\end{equation}
where $id$ denotes the identity operator on $\mathcal{A}(\Omega)$.

\begin{theorem}\label{tULleqLU} For any $\delta>0$ we have
$U_\delta\circ L_\delta\leq L_\delta\circ U_\delta$.
\end{theorem}
\begin{proof} Let $f\in\mathcal{A}(\Omega)$ and let $x\in\Omega$. Denote
$p=(L_\delta\circ
U_\delta)(f)(x)=S_\frac{\delta}{2}(I_\delta(S_\frac{\delta}{2}(f)))(x)$.
Let $\varepsilon$ be an arbitrary positive. For every $y\in
B_\frac{\delta}{2}(x)$ we have
\begin{equation}\label{ineqtULleqLU}
I_\delta(S_\frac{\delta}{2}(f))(y)\leq p<p+\varepsilon.
\end{equation}
\underline{Case 1.} There exists $z\in B_\frac{\delta}{2}(x)$ such
that $S_\frac{\delta}{2}(f)(z)<p+\varepsilon$. Then
$f(t)<p+\varepsilon \ \mbox{ for }\ t\in B_\frac{\delta}{2}(z)$,
which implies that $I_\frac{\delta}{2}(f)(t)< p+\varepsilon \
\mbox{ for }\ t\in B_\delta(z)$. Hence
$S_\delta(I_\frac{\delta}{2}(f))(z)\leq p+\varepsilon$. Then
$(U_\delta\circ
L_\delta)(f)(t)=I_\frac{\delta}{2}(S_\delta(I_\frac{\delta}{2}(f)))(t)\leq
p+\varepsilon$ for $t\in B_\frac{\delta}{2}(z)$. Since $x\in
B_\frac{\delta}{2}(z)$, see the case assumption, from the above
inequality we have $(U_\delta\circ L_\delta)(f)(x)\leq
p+\varepsilon$.\\
\underline{Case 2.} For every $z\in
B_\frac{\delta}{2}(x)$ we have $S_\frac{\delta}{2}(f)(z)\geq
p+\varepsilon$. Denote
$D=\left\{z\in\Omega:S_\frac{\delta}{2}(f)(z)< p+\varepsilon
\right\}$. We will show that for every $z\in B_\delta(x)$ we have
\begin{equation}\label{intersecttULleqLU}
B_\delta(z)\cap D\neq\emptyset
\end{equation}
Due to the inequality (\ref{ineqtULleqLU}) we have that
(\ref{intersecttULleqLU}) holds for every $z\in
B_\frac{\delta}{2}(x)$. Let $z\in B_\delta(x)$ and let
$z>x+\frac{\delta}{2}$. This implies that
$x+\frac{\delta}{2}\in\Omega$. Using the inequality
(\ref{ineqtULleqLU}) for $y=x+\frac{\delta}{2}$ as well as the
case assumption we obtain that the set
$\left(x+\frac{\delta}{2},x+\frac{3\delta}{2}\right]\cap D$ is not
empty. Then $B_\delta(z)\cap D\supset
\left(x+\frac{\delta}{2},x+\frac{3\delta}{2}\right]\cap D\neq
\emptyset$. For $z<x-\frac{\delta}{2}$ condition
(\ref{intersecttULleqLU}) is proved in a similar way. Hence
(\ref{intersecttULleqLU}) holds for all $z\in B_\delta(x)$. Let
$z\in B_\delta(x)$ and $v\in B_\delta(y)\cap D$. Since $v\in D$ we
have $f(t)< p+\varepsilon$, for $t\in B_\frac{\delta}{2}(v)$.
Using that $B_\frac{\delta}{2}(z)\cap B_\frac{\delta}{2}(v)\neq
\emptyset$ we obtain that $I_\frac{\delta}{2}(f)(z)<
p+\varepsilon$, $z\in B_\delta(x)$. Therefore
$S_\delta(I_\frac{\delta}{2}(f))(x)\leq p+\varepsilon$. Then
\[
(U_\delta\circ
L_\delta)(f)(x)=I_\frac{\delta}{2}(S_\delta(I_\frac{\delta}{2}(f)))(x)\leq
S_\delta(I_\frac{\delta}{2}(f))(x)\leq p+\varepsilon.
\]
Combining the results of Case 1 and Case 2 we have $(U_\delta\circ
L_\delta)(f)(x)\leq p+\varepsilon$. Since $\varepsilon$ is
arbitrary this implies that $(U_\delta\circ L_\delta)(f)(x)\leq
p=(L_\delta\circ U_\delta)(f)(x)$.
\end{proof}

\begin{theorem}\label{tLUULidemp}
For a given $\delta>0$ the operators $L_\delta\circ U_\delta$ and
$U_\delta\circ L_\delta$ are both idempotent.
\end{theorem}
The proof is an immediate application of Lemma 1.

\begin{theorem}\label{tULU=LU}
We have $U_\delta\circ L_\delta\circ U_\delta=L_\delta\circ
U_\delta$, $L_\delta\circ U_\delta\circ L_\delta=U_\delta\circ
L_\delta$, $\delta>0$.
\end{theorem}
\begin{proof}
Using the inequalities (\ref{LleqidleqU}) and the monotonicity of
the operators $L_\delta$, $U_\delta$, see (\ref{LdUdfmon}), we
obtain $U_\delta\circ L_\delta\circ U_\delta \geq id\circ
L_\delta\circ U_\delta\ =\ L_\delta\circ U_\delta$. For the proof
of the inverse inequality we use Theorem \ref{tULleqLU} and the
idempotence of $U_\delta$ as follows:
\[
U_\delta\circ L_\delta\circ U_\delta=(U_\delta\circ L_\delta)\circ
U_\delta\leq (L_\delta\circ U_\delta)\circ U_\delta=L_\delta\circ
(U_\delta\circ U_\delta)=L_\delta\circ U_\delta
\]
Therefore $U_\delta\circ L_\delta\circ U_\delta=L_\delta\circ
U_\delta$. The second equality is proved in a similar way.
\end{proof}

It follows from Theorems \ref{tLUULidemp} and \ref{tULU=LU} that
for a fixed $\delta>0$ every composition involving finite number
of the operators $L_\delta$ and $U_\delta$ is an element of the
set $\{L_\delta,U_\delta,U_\delta\circ L_\delta, L_\delta\circ
U_\delta\}$. Hence the operators $L_\delta$ and $U_\delta$ form a
semi-group with a composition table as follows:
\begin{center}
\begin{tabular}{c||c|c|c|c|}
&$ L_\delta$&$ U_\delta$&$U_\delta\circ L_\delta$&$L_\delta\circ
U_\delta$\\\hline\hline $L_\delta$&$ L_\delta$&$
 L_\delta\circ U_\delta$&$U_\delta\circ L_\delta$
&$ L_\delta\circ U_\delta$\\\hline $ U_\delta$&$ U_\delta\circ
L_\delta$&$ U_\delta$&$ U_\delta\circ L_\delta$ &$L_\delta\circ
U_\delta$\\\hline $ U_\delta\circ L_\delta$&$ U_\delta\circ
L_\delta$&$L_\delta\circ U_\delta$&$U_\delta\circ
L_\delta$&$L_\delta\circ U_\delta$\\\hline $ L_\delta\circ
U_\delta$&$U_\delta\circ L_\delta$&$ L_\delta\circ U_\delta$
&$U_\delta\circ L_\delta$ &$ L_\delta\circ U_\delta$\\\hline
\end{tabular}
\end{center}
Furthermore, an easy application of Theorem \ref{tULleqLU} shows
that this semi-group is completely ordered. Namely, we have
$L_\delta\leq U_\delta\circ L_\delta\leq L_\delta\circ
U_\delta\leq U_\delta$.

\section{Locally $\delta$-monotone approximations}

\begin{definition}
Let $\delta>0$ and a function $f\in\mathbb{A}(\Omega)$ be given.

(i) The function $f$ is called upwards $\delta$-monotone if for
every interval $[x,y]\subset\Omega$ with $y-x\leq\delta$ we have
$\displaystyle \sup_{z\in[x,y]} f(z)=\max\{f(x),f(y)\}$.

(ii) The function $f$ is called downwards $\delta$-monotone if for
every interval $[x,y]\subset\Omega$ with $y-x\leq\delta$ we have
$\displaystyle \inf_{z\in[x,y]} f(z)=\min\{f(x),f(y)\}$.

(iii) The function $f$ is called locally $\delta$-monotone if it
is both downwards $\delta$-monotone and upwards $\delta$-monotone.
\end{definition}

The name locally $\delta$-monotone reflects the following
characterization:
\begin{equation}\label{dmonchar}
\begin{tabular}{l}$f$ is locally\\ $\delta$-monotone\end{tabular}
\!\Longleftrightarrow\!
\begin{tabular}{l}On any interval $[x,y]\!\subseteq\! \Omega$, $y\!-\!x\!\leq\!\delta$, $f$ is either\\
monotone increasing or monotone decreasing
\end{tabular}
\end{equation}

\begin{theorem}
For every $\delta>0$ and $f\in\mathcal{A}(\Omega)$ the function
$S_\frac{\delta}{2}(f)$ is upwards $\delta$-monotone while the
function $I_\frac{\delta}{2}(f)$ is downwards $\delta$-monotone.
\end{theorem}
\begin{proof}
Let $\delta>0$, $f\in\mathcal{A}(\Omega)$ and
$[x,y]\subseteq\Omega$, $y-x\leq\delta$. Denote
$g=S_\frac{\delta}{2}(f)$. It is easy to see that for every
$z\in[x,y]$ we have $B_\frac{\delta}{2}(z)\subseteq
B_\frac{\delta}{2}(x)\cup B_\frac{\delta}{2}(y)$. Therefore
\begin{eqnarray*}
g(z)&=&\sup\{f(t):t\in B_\frac{\delta}{2}(z)\}\ \leq\
\sup\{f(t):t\in B_\frac{\delta}{2}(x)\cup
B_\frac{\delta}{2}(y)\}\\
&=&\max\{\sup\{f(t):t\in B_\frac{\delta}{2}(x)\},\sup\{f(t):t\in
B_\frac{\delta}{2}(z)\}\}\\
&=&\max\{g(x),g(y)\},
\end{eqnarray*}
which shows that function $g$ is upwards $\delta$-monotone. The
downwards $\delta$-monotonicity of $U_\delta(f)$ is proved in a
similar way.
\end{proof}

\begin{corollary}\label{tLdupmon}
For every $\delta>0$ and $f\in\mathcal{A}(\Omega)$ the function
$L_\delta(f)$ is upwards $\delta$-monotone while the function
$U_\delta(f)$ is downwards $\delta$-monotone.
\end{corollary}

\begin{theorem}\label{tULlocdmon}
For every $\delta>0$ and $f\in\mathcal{A}(\Omega)$ the functions
$(U_\delta\circ L_\delta)(f)$ and $(L_\delta\circ U_\delta)(f)$
are both locally $\delta$-monotone.
\end{theorem}
The proof follows from Corollary \ref{tLdupmon} and the
composition table in the preceding section. Theorem
\ref{tULlocdmon} shows that the operators $U_\delta\circ L_\delta$
and $L_\delta\circ U_\delta$ provide locally $\delta$-monotone
approximations to the functions in $\mathcal{A}(\Omega)$. The
error of the approximation can be estimated in terms of the
modulus of nonmonotonicity. Let us recall the definition.

\begin{definition}\label{defmodulus}
Let $f\!\in\!\mathcal{A}(\Omega)$. The mapping
$\mu(f,\cdot):\mathbb{R}^+\!\rightarrow\!
\mathbb{R}^+\!\cup\!\{0\}$ given by
\[
\mu(f,\delta)\!=\!\frac{1}{2}\!\!\!\!
\sup_{\begin{tabular}{c}$x_{1,2}\in\Omega$\\$0\!<\!x_2\!-\!x_2\!\leq\!\delta$\end{tabular}}
\!\!\!\!\!\sup_{\begin{tabular}{c}$x\in[x_1,x_2]$\end{tabular}}
\!\!(|f(x_1)\!-\!f(x)|\!+\!|f(x_2)\!-\!f(x)|\!-\!|f(x_1)\!-\!f(x_2)|)
\]
is called modulus of nonmonotonicity of $f$.
\end{definition}
The locally $\delta$-monotone functions can be conveniently
characterized through the modulus of nonmonotonicity. For any
$f\in\mathcal{A}(\Omega)$ and $\delta>0$ we have
\begin{equation}\label{dmonchar2}
f \mbox{ is locally }\delta\mbox{-monotone}\ \Longleftrightarrow \
\mu(f,\delta)=0
\end{equation}

We will derive error estimates first in the case when
$\Omega=\mathbb{R}$. It will prove useful to consider the upper
semi-continuous envelope of the modulus of nonmonotonicity. Let
$f\in\mathcal{A}(\Omega)$. Using that $\mu(f,\delta)$ is monotone
increasing with respect to $\delta$ the upper semi-continuous
envelope of $\mu$ can be represented as
$\hat{\mu}(f,\delta)=S(\mu(f,\cdot))(\delta)=\lim_{\varepsilon\rightarrow
0^+}\mu(f,\delta+\varepsilon)$.
\begin{theorem}\label{testim1}
Let $f\in\mathcal{A}(\mathbb{R})$ and $\delta>0$. Then
$f(x)-L_\delta(f)(x)\leq\hat\mu(f,\delta)$,
$U_\delta(f)(x)-f(x)\leq\hat\mu(f,\delta)$, $x\in\Omega$.
\end{theorem}
\begin{proof}
Let $x\in\mathbb{R}$. Denote $p=L_\delta(f)(x)$. If $p=f(x)$ the
first inequality of the theorem holds. Assume that $p<f(x)$. Let
$\eta>0$ be such that $p+\eta<f(x)$. Then we have
\begin{equation}\label{ineq1testim1}
I_\frac{\delta}{2}(f)(y)<p+\eta,\ y\in B_\frac{\delta}{2}(x).
\end{equation}
Denote $D_1=\{z\leq x:f(z)<p+\eta\}$, $D_2=\{z\geq
x:f(z)<p+\eta\}$, $z_1=\sup D_1$, $z_2=\sup D_2$. Using the
inequality (\ref{ineq1testim1}) with $y=x-\frac{\delta}{2}$ and
$y=x+\frac{\delta}{2}$ we obtain that $D_1$ and respectively $D_2$
are not empty and that $x-\delta\leq z_1\leq x\leq z_2\leq
x+\delta$. Therefore
$z_3=\frac{z_1+z_2}{2}\in\left[\frac{x-\delta+x}{2},\frac{x+x+\delta}{2}\right]
=B_\frac{\delta}{2}(x)$.
Then the inequality (\ref{ineq1testim1}) implies that
$B_\frac{\delta}{2}\cap (D_1\cup D_2)\neq\emptyset$. Hence
$z_2-z_3=z_3-z_1\leq\frac{\delta}{2}$.

Let $\varepsilon>0$ be arbitrary. The neighborhood
$B_\frac{\delta+\varepsilon}{2}(z_3)$ has nonempty intersections
with both $D_1$ and $D_2$. Let $t_1\in
B_\frac{\delta+\varepsilon}{2}\cap D_1$ and $t_2\in
B_\frac{\delta+\varepsilon}{2}\cap D_2$. We have
$t_2-t_1<\delta+\varepsilon$ and $x\in[t_1,t_2]$. From the
definition of the modulus of nonmonotonicity we have
$|f(t_1)-f(x)|+|f(t_2)-f(x)|-|f(t_1)-f(t_2)|\leq
2\mu(f,\delta+\varepsilon)$. On the other side
$|f(t_1)-f(x)|+|f(t_2)-f(x)|-|f(t_1)-f(t_2)|=2f(x)-f(t_1)-f(t_2)-|f(t_1)-f(t_2)|
=2f(x)-\max\{f(t_1),f(t_2)\}>2f(x)-2(p+\eta)$. Therefore
$f(x)-p-\eta<\mu(f,\delta+\varepsilon)$. Going with $\varepsilon$
to $0$ we obtain $f(x)-p\leq\hat\mu(f,\delta)+\eta$. Since $\eta$
is arbitrary small this implies the first inequality of the
Theorem. The second inequality is proved in a similar manner.
\end{proof}

Using Theorem \ref{testim1} as well as (\ref{dmonchar}) and
Corollary \ref{tLdupmon} we have the following characterization of
the fixed points of operators $L_\delta$ and $U_\delta$. For any
$f\in\mathcal{A}(\mathbb{R})$
\[
\hat\mu(f,\delta)=0\ \Longrightarrow\ (L_\delta(f)=f,\
U_\delta(f)=f)\ \Longrightarrow\ \mu(f,\delta)=0
\]

\begin{theorem}\label{testim2}
Let $\delta>0$ and $f\in\mathcal{A}(\mathbb{R})$ . Then
$\mu(L_\delta(f),\delta)\leq\mu(f,\delta)$,
$\hat\mu(L_\delta(f),\delta)\leq\hat\mu(f,\delta)$,
$\mu(U_\delta(f),\delta)\leq\mu(f,\delta)$,
$\hat\mu(U_\delta(f),\delta)\leq\hat\mu(f,\delta)$.
\end{theorem}
\begin{proof}
We will prove the inequalities for $L_\delta(f)$ since the ones
for $U_\delta(f)$ are proved in a similar way. Denote
$g=L_\delta(f)$. Let $[x_1,x_2]$ be an arbitrary interval of
length at most $\delta$ and let $x\in[x_1,x_2]$. We consider the
number $q=|f(x_1)-f(x)|+|f(x_2)-f(x)|-|f(x_1)-f(x_2)|$. According
to Corollary \ref{tLdupmon} the function $g$ is upper $\delta$
monotone, which implies that $g(x)\leq\max\{g(x_1),g(x_2)\}$. If
we also have $g(x)\geq\min\{g(x_1),g(x_2)\}$, then the number $q$
is zero and the first inequality of the Theorem is trivially
satisfies. Let $g(x)<\min\{g(x_1),g(x_2)\}$ and let $\eta>0$ be
such that $g(x)+\eta<\min\{g(x_1),g(x_2)\}$. The number $q$ can
then be represented in the form $q=2(\min\{g(x_1),g(x_2)\}-g(x))$.
If we assume that $f(y)\geq g(x)+\eta$ for all $y\in
B_\frac{\delta}{2}(x)$, then $g(x)\geq
I_\frac{\delta}{2}(f)(x)\geq g(x)+\eta$, which is a contradiction.
Therefore, there exists $y\in B_\frac{\delta}{2}(x)$ such that
$f(y)<g(x)+\eta$. If $y<x_1$ then using that
$B_\frac{\delta}{2}(x_1)\subseteq B_\frac{\delta}{2}(y)\cup
B_\frac{\delta}{2}(x)$ we obtain $I_\frac{\delta}{2}(z)<g(x)+\eta$
for all $z\in B_\frac{\delta}{2}(x_1)$ which implies
$g(x_1)=S_\frac{\delta}{2}(I_\frac{\delta}{2}(f))(x_1)\leq
g(x)+\eta<g(x_1)$. This contradiction shows that $y\geq x_1$. In a
similar way we show that $y\leq x_2$. Using also the first
inequality of Theorem \ref{tLdfleqf} we have
$q=2(\min\{g(x_1),g(x_2)\}-g(x))\leq
2(\min\{f(x_1),f(x_2)\}-f(y)-2\eta)\leq\mu(f,\delta)-2\eta$. Using
that $\eta$ can be arbitrary small we obtain $q\leq
\mu(f,\delta)$. Since the interval $[x_1,x_2]$ of length at most
$\delta$ and $x\in[x_1,x_2]$ are arbitrary this implies that
$\mu(g,\delta)\leq\mu(f,\delta)$. The inequality
$\hat\mu(g,\delta)\leq\hat\mu(f,\delta)$ follows from the
monotonicity of the operator $S$, see (\ref{IdSdmon}).
\end{proof}

In the next theorem we give estimates for the error of
approximation of a function $f\in\mathcal{A}(\mathbb{R})$ in terms
of the supremum norm denoted here by $||\cdot||$.

\begin{theorem}\label{testim3}
Let $\delta>0$ and $f\in\mathcal{A}(\mathbb{R})$ . Then
$||f-(L_\delta\circ U_\delta)(f)||\leq \hat\mu(f,\delta)$,
$||f-(U_\delta\circ L_\delta)(f)||\leq \hat\mu(f,\delta)$.
\end{theorem}
\begin{proof}
Applying Theorems \ref{testim1} and \ref{testim2} we obtain
$(L_\delta\circ U_\delta)(f)(x)\leq U_\delta(f)(x)\ \leq\
f(x)+\hat\mu(f,\delta)$ and $(L_\delta\circ U_\delta)(f)(x)\geq
U_\delta(f)(x)-\hat\mu(L_\delta(f),\delta)\ \geq\
f(x)-\hat\mu(f,\delta)$, which implies the first inequality of the
Theorem. The second inequality is proved in a similar way.
\end{proof}

Error estimates similar to Theorem \ref{testim3} can be derived in
case of $\Omega$ being finite or semi-finite interval using a
modification of the modulus of nonmonotonicity. For simplicity we
will only consider the case $\Omega=[a,b]$, $a,b\in\mathbb{R}$.

\begin{definition}
Let $f\!\in\!\mathcal{A}(\Omega)$. The mapping
$\tilde\mu(f,\cdot):\mathbb{R}^+\!\!\rightarrow\!\!
\mathbb{R}^+\!\cup\!\{0\}$ given by
\[
\tilde\mu(f,\delta)=\sup\left\{\hat\mu(f,\delta),
\sup_{x_{1,2}\in[a,a+\frac{\delta}{2}]}(|f(x_1)\!-\!f(x_2)|,
\sup_{x_{1,2}\in[b-\frac{\delta}{2},b]}(|f(x_1)\!-\!f(x_2)|\right\}
\]
is called modified modulus of nonmonotonicity of $f$.
\end{definition}

This modulus is similar to the corrected modulus of nonmotonicity
in \cite{Markov}.

\begin{theorem}\label{testim4}
Let $f\in\mathcal{A}(\mathbb{R})$ and $\delta>0$. Then
$f(x)-L_\delta(f)(x)\leq\tilde\mu(f,\delta)$,
$U_\delta(f)(x)-f(x)\leq\tilde\mu(f,\delta)$, $x\in\Omega$.
\end{theorem}
The proof is similar to the proof of Theorem \ref{testim1}.

It is easy to see that for any $f\in\mathcal{A}[a,b]$ the
functions $L_\delta(f)$ and $U_\delta(f)$ are constants on each of
the intervals $[a,a+\frac{\delta}{2}]$ and $[b-\frac{\delta}{2}]$.
Therefore, using also Theorem \ref{testim2} we have
$\tilde{\mu}(L_\delta(f),\delta)=\hat\mu(L_\delta(f),\delta)\leq
\hat\mu(f,\delta)\leq\tilde\mu(f,\delta)$. In the same way we
obtain $\tilde{\mu}(U_\delta(f),\delta)\leq\tilde\mu(f,\delta)$.
Hence the modified modulus satisfies the similar inequalities to
the ones given in Theorem \ref{testim2} for $\mu$ and $\hat\mu$.
Using these inequalities and Theorem \ref{testim4} we obtain the
error estimates in the next theorem, which are similar to the ones
in Theorem \ref{testim3}.

\begin{theorem}\label{testim5}
Let $\delta>0$ and $f\in\mathcal{A}(\mathbb{R})$ . Then
$||f-(L_\delta\circ U_\delta)(f)||\leq \tilde\mu(f,\delta)$,
$||f-(U_\delta\circ L_\delta)(f)||\leq \tilde\mu(f,\delta)$.
\end{theorem}

\section{Conclusion}
In this paper we extended the LULU operators from sequences to
real functions defined on a real interval using the lower and
upper $\delta$-envelopes of functions. The obtained structure,
although more general than the well known LULU structure of the
discrete operators, retains some of its essential properties. For
a fixed $\delta>0$ the compositions $L_\delta\circ U_\delta$ and
$U_\delta\circ L_\delta$ provide locally $\delta$-monotone
approximations for real functions, the error of approximation
being estimated in terms of the modulus of nonmonotonicity of the
functions. Further properties of the LULU operators for functions
on continuous domains, e.g. trend preservation, will be
investigated in the future. Generalizing the theory to functions
on multidimensional domains is still an open problem.


\begin{thebibliography}{9}

\bibitem{Markov} S. Markov, Relations between the integral and
Hausdorff distance with applications to differential equations,
Pliska {\bf 1} (1977) 112--121.

\bibitem{RohwerQM2002} C. H. Rohwer, Variation reduction and
$LULU$-smoothing, Quaestiones Mathematicae {\bf 25} (2002)
163--176.

\bibitem{RohwerQM2004} C. H. Rohwer, Fully trend preserving operators,
Quaestiones Mathematicae {\bf 27} (2004) 217--230.


\bibitem{Rohwerbook} C. H. Rohwer, Nonlinear Multiresolution
Analysis, Birkh\"{a}user, 2005.

\bibitem{Sendov} B. Sendov, Hausdorff Approximations, Kluwer,
Boston, 1990.

\bibitem{Serra} J. Serra, Image Analysis and Mathematical
Morphology, Academic Press, London, 1982.

\end{thebibliography}
\end{document}